\documentclass[a4paper,12pt]{article}
\usepackage{amsmath,amsfonts,amssymb}
\usepackage[all]{xy}
\makeatletter
\newbox\bk@bxb
\newbox\bk@bxa
\newif\if@bkcont
\newcount\bk@lcnt

\def\breakboxskip{2pt}
\def\breakboxparindent{1.8em}

\def\breakbox{\vskip\breakboxskip\relax
\setbox\bk@bxb\vbox\bgroup
\advance\linewidth -2\fboxrule
\hsize\linewidth\@parboxrestore
\parindent\breakboxparindent\relax}

\def\bk@split{%
\@tempdimb\ht\bk@bxb 
\advance\@tempdimb\dp\bk@bxb
\setbox\bk@bxa\vsplit\bk@bxb to\z@ 
\setbox\bk@bxa\vbox{\unvbox\bk@bxa}
\setbox\@tempboxa\vbox{\copy\bk@bxa\copy\bk@bxb}
\advance\@tempdimb-\ht\@tempboxa
\advance\@tempdimb-\dp\@tempboxa}

\def\bk@addfsepht{%
\setbox\bk@bxa\vbox{\vskip\fboxsep\box\bk@bxa}}

\def\bk@addskipht{%
\setbox\bk@bxa\vbox{\vskip\@tempdimb\box\bk@bxa}}

\def\bk@addfsepdp{%
\@tempdima\dp\bk@bxa
\advance\@tempdima\fboxsep
\dp\bk@bxa\@tempdima}

\def\bk@addskipdp{%
\@tempdima\dp\bk@bxa
\advance\@tempdima\@tempdimb
\dp\bk@bxa\@tempdima}

\def\bk@line{%
\hbox to \linewidth{%
\hskip-2\fboxsep\vrule \@width\fboxrule\hskip.5\fboxsep\vrule \@width\fboxrule\hskip1.5\fboxsep
\box\bk@bxa\hfil
}}%

\def\endbreakbox{\egroup
\ifhmode\par\fi{\noindent\bk@lcnt\@ne
\@bkconttrue\baselineskip\z@\lineskiplimit\z@
\lineskip\z@\vfuzz\maxdimen
\bk@split\bk@addfsepht\bk@addskipdp
\ifvoid\bk@bxb 
\def\bk@fstln{\bk@addfsepdp
\hskip-\parindent\vbox{\llap{\raisebox{-2ex}{\rule{1.5\fboxsep}{\fboxrule}\hskip.5\fboxsep}}\bk@line\llap{\rule{1.5\fboxsep}{\fboxrule}\hskip.5\fboxsep}}}

\else 
\def\bk@fstln{\vbox{\llap{\raisebox{-2ex}{\rule{1.5\fboxsep}{\fboxrule}\hskip.5\fboxsep}}\bk@line}\hfil%
\advance\bk@lcnt\@ne
\loop
\bk@split\bk@addskipdp\leavevmode
\ifvoid\bk@bxb 
\@bkcontfalse\bk@addfsepdp
\vtop{\bk@line\llap{\rule{2\fboxsep}{\fboxrule}}}%

\else 
\bk@line
\fi
\hfil\advance\bk@lcnt\@ne
\if@bkcont\repeat}%
\fi
\leavevmode\bk@fstln\par}\vskip\breakboxskip\relax}

\def\smp{\smallskip\par}

\def\pf{\noindent{\bf Proof~:}\ }

\def\findemo{~\leaders\hbox to 1em{\hss\  \hss}\hfill~\raisebox{.5ex}{\framebox[1ex]{}}\smp}

\def\mpn{\medskip\par\noindent}
\def\smpn{\smallskip\par\noindent}

\def\smp{\smallskip\par}
\def\smpn{\smallskip\par\noindent}

\def\mpoint{\;\;.}
\def\mvirg{\;\;,}

\def\Res{{\rm Res}}
\def\Ind{{\rm Ind}}

\def\Hom{{\rm Hom}}
\def\End{{\rm End}}

\def\Im{{\rm Im}}

\def\Id{{\rm Id}}

\def\Z{\mathbb{Z}}

\newcommand{\romain}[1]{\uppercase\expandafter{\romannumeral #1}}

\newcommand{\flh}[2]{\mathop{\hbox to 12mm{\rightarrowfill}}_{\displaystyle #2}^{\displaystyle #1}\limits}
\newcommand{\sflh}[2]{\mathop{\hbox to 12mm{\rightarrowfill}}_{\scriptstyle #2}^{\scriptstyle #1}\limits}

\newcommand{\gsetover}[2]{#1\hbox{-{\bf set}}\smash{\downarrow}_{#2}}

\newcommand{\gset}[1]{#1\hbox{-$\mathsf{set}$}}

\newcommand{\gMod}[1]{#1{\hbox{-}\mathsf{Mod}}}

\newcommand{\sou}[1]{\underline{#1}}

\newcommand{\sumb}[2]{\mathop{\sum}_{{\scriptstyle #1}\atop {\scriptstyle #2}}}

\newcommand{\carre}[8]{\begin{array}{ccc}
#1&\mathop{\hbox to 12mm{\rightarrowfill}}^{\displaystyle{#2}}\limits&#3\\
\llap{$\displaystyle{#4}$}\left\downarrow\vbox to 6mm{}\right. & & \left\downarrow\vbox to 6mm{}\right.\rlap{$\displaystyle{#5}$}\\
#6&\mathop{\hbox to 12mm{\rightarrowfill}}_{\displaystyle #7}\limits&#8\\
\end{array}}

\newcommand{\carrem}[8]{\begin{array}{ccc}
#1&\mathop{\hbox to 12mm{\rightarrowfill}}^{\displaystyle #2}\limits&#3\\
\llap{$\displaystyle #4$}\left\uparrow\vbox to 6mm{}\right. & & \left\uparrow\vbox to 6mm{}\right.\rlap{$\displaystyle #5$}\\
#6&\mathop{\hbox to 12mm{\rightarrowfill}}_{\displaystyle #7}\limits&#8\\
\end{array}}

\newenvironment{enonce}[1]{\pagebreak[2]\refstepcounter{subsection}\refstepcounter{prop}\smpn{{\bf \thesection.\arabic{prop}.\ \ #1~:}}\begin{it} }{\end{it}\smp}
\newenvironment{enonce*}[1]{\pagebreak[2]\smpn{#1~:}\begin{it} }{\end{it}\smp}
\newcommand{\result}[1]{\begin{enonce}{#1}}
\def\fresult{\end{enonce}}
\newcommand{\npar}{\smallskip\par\noindent\pagebreak[2]\refstepcounter{subsection}\refstepcounter{prop}{\bf \thesection.\arabic{prop}.\ \ }}
\newcommand{\masubsect}[1]{\medskip\par\noindent\pagebreak[3]\refstepcounter{subsection}\refstepcounter{prop}{\bf \thesection.\arabic{prop}.\ \ #1.\ }}

\newenvironment{mth}[1]{\begin{breakbox}\begin{enonce}{#1}}{\end{enonce}\end{breakbox}}
\newenvironment{mth*}[1]{\begin{breakbox}\begin{enonce*}{#1}}{\end{enonce*}\end{breakbox}}
\newenvironment{rem}[1]{\refstepcounter{subsection}\refstepcounter{prop} \mpn{{\bf \thesection.\arabic{prop}.}\ \ \bf#1\ :}}{\smp}

\def\dom{\backslash}
\makeatletter
\renewenvironment{enumerate}{\ifnum \@enumdepth >3 \@toodeep\else
      \advance\@enumdepth \@ne
      \edef\@enumctr{enum\romannumeral\the\@enumdepth}\list
      {\csname label\@enumctr\endcsname}{\setlength{\topsep}{1ex}\setlength{\itemsep}{0pt}\usecounter
        {\@enumctr}\def\makelabel##1{\hss\llap{##1}}}\fi}{\endlist}
\renewenvironment{itemize}{\ifnum \@itemdepth >3 \@toodeep\else \advance\@itemdepth \@ne
\edef\@itemitem{labelitem\romannumeral\the\@itemdepth}%
\list{\csname\@itemitem\endcsname}{\setlength{\topsep}{1ex}\setlength{\itemsep}{0pt}\def\makelabel##1{\hss\llap{##1}}}\fi}
{\endlist}
\makeatother
\makeatletter
\def\@sect#1#2#3#4#5#6[#7]#8{\ifnum #2>\c@secnumdepth
    \let\@svsec\@empty\else
    \refstepcounter{#1}\edef\@svsec{\csname the#1\endcsname .\hskip .5em}\fi
    \@tempskipa #5\relax
     \ifdim \@tempskipa>\z@
       \begingroup #6\relax
         \@hangfrom{\hskip #3\relax\@svsec}{\interlinepenalty \@M #8\par}%
       \endgroup
      \csname #1mark\endcsname{#7}\addcontentsline
        {toc}{#1}{\ifnum #2>\c@secnumdepth \else
                     \protect\numberline{\csname the#1\endcsname}\fi
                   #7}\else
       \def\@svsechd{#6\hskip #3\relax  
                  \@svsec #8\csname #1mark\endcsname
                     {#7}\addcontentsline
                          {toc}{#1}{\ifnum #2>\c@secnumdepth \else
                            \protect\numberline{\csname the#1\endcsname}\fi
                      #7}}\fi
    \@xsect{#5}}
\def\section{\@startsection {section}{1}{\z@}{-3.5ex plus-1ex minus
    -.2ex}{2.3ex plus.2ex}{\reset@font\Large\bf}}  

\makeatother
\renewenvironment{equation}{\refstepcounter{subsection}\refstepcounter{prop}$$}{\leqno{\bf (\theprop)}$$}

\def\mar[#1]{\ar@{-}[#1]|-{\object@{<}}}
\def\marb[#1]{\ar@{-}[#1]|{\object+{  }}}

\newcommand{\sgset}[1]{#1\hbox{-\underline{{\bf set}}}}
\newcommand{\bisetover}[3]{(#1,#2)\hbox{-{\bf biset}}\smash{\downarrow}_{#3}}
\newcommand{\igset}[1]{#1\hbox{-}\sf{set}}
\newcommand{\isgset}[1]{#1\hbox{-}\sou{\sf{set}}}

\begin{document}
\centerline{\LARGE\bf Fused Mackey functors}\vspace{.5cm}
\centerline{Serge Bouc}\vspace{.5cm}
{\footnotesize {\bf Abstract:} Let $G$ be a finite group. In \cite{htw_Mackey-bisets}, Hambleton, Taylor and Williams have considered the question of comparing Mackey functors for $G$ and biset functors defined on subgroups of $G$ and bifree bisets as morphisms.\par }
{\footnotesize This paper proposes a different approach to this problem, from the point of view of various categories of $G$-sets. In particular, the category $\sgset{G}$ of {\em fused $G$-sets} is introduced, as well as the category $\sou{\mathbf{S}}(G)$ of spans in $\sgset{G}$. The {\em fused Mackey functors} for $G$ over a commutative ring $R$ are defined as $R$-linear functors from $R\,\sou{\mathbf{S}}(G)$ to $R$-modules. They form an abelian subcategory $\mathsf{Mack}_R^f(G)$ of the category of Mackey functors for $G$ over~$R$. The category $\mathsf{Mack}_\Z^f(G)$ is equivalent to the category of conjugation Mackey functors of~\cite{htw_Mackey-bisets}. The category $\mathsf{Mack}_R^f(G)$ is also equivalent to the category of modules over the {\em fused Mackey algebra} $\mu_R^f(G)$, which is a quotient of the usual Mackey algebra $\mu_R(G)$ of $G$ over $R$.\vspace{1ex}\par}
{\footnotesize {\bf AMS Subject classification :} 18A25, 19A22, 20J15.\vspace{1ex}}\par

\section{Introduction}
This note is devoted to the frequently asked question of comparing Mackey functors for a single finite group $G$ with biset functors defined only on subgroups of $G$ and left-right free bisets as morphisms. The answer to this question has already been given by Hambleton, Taylor and Williams (\cite{htw_Mackey-bisets}), but in a rather computational and non canonical way (in particular, in Section 7, the definition of the functor $j_\bullet$ requires the choice of sets of representatives of orbits of any finite $G$-set). \par
The present paper makes a systematic use of Dress definition (\cite{dress}) and Lindner definition (\cite{lindner}) of Mackey functors, to avoid these non canonical choices. This leads to the definition of the category of {\em fused $G$-sets} (Section~\ref{fused G-sets}), and the category of {\em fused Mackey functors} (Section~\ref{fused Mackey functors}) for a finite group $G$, which is equivalent to the category of ``conjugation invariant Mackey functors" of~\cite{htw_Mackey-bisets}. This category is also equivalent to the category of modules over the {\em fused Mackey algebra}, introduced in Section~\ref{algebra}.
\par
\section{Conjugation bisets revisited}
\npar First a notation~: when $G$ is a finite group, and $X$ is a finite $G$-set, let $\gsetover{G}{X}$ denote the category of (finite) $G$-sets over $X$: its objects are pairs $(Y,b)$ consisting of a finite $G$-set $Y$, and a morphism of $G$-sets $b:Y\to X$. A morphism $f:(Y,b)\to (Z,c)$ in $\gsetover{G}{X}$ is a morphism of $G$-sets $f:Y\to Z$ such that $c\circ f=b$.\par
There is an obvious notion of disjoint union in $\gsetover{G}{X}$, and the corresponding Grothendieck group is called the Burnside group over $X$. It will be denoted by $\mathcal{B}(_GX)$, or $\mathcal{B}(X)$ when~$G$ is clear from the context.\par
Similarly, when $G$ and $H$ are finite groups, and $U$ is a $(G,H)$-biset, one can define the category $\bisetover{G}{H}{U}$ of $(G,H)$-bisets over $U$, and the Burnside group $\mathcal{B}(_GU_H)$ of $(G,H)$-bisets over $U$.
\npar \label{ind}When $H$ is a subgroup of $G$, and $Y$ is an $H$-set, induction from $H$-sets to $G$-sets is an equivalence of categories from $\gsetover{H}{Y}$ to $\gsetover{G}{\Ind_H^GY}$. A quasi-inverse equivalence is the functor 
sending the $G$-set $(X,a)$ over $\Ind_H^GY$ to the $H$-set $a^{-1}(1\times_HY)$ (see \cite{green}~Lemma 2.4.1). In particular $\mathcal{B}(_HY)\cong\mathcal{B}(_G\Ind_H^GY)$.
\npar \label{observe}Now an observation: when $H$ and $K$ are subgroups of $G$, the conjugation $(K,H)$-bisets defined in Section~6 of \cite{htw_Mackey-bisets} are exactly those over the biset $_KG_H$ (the set $G$ on which $K$ and $H$ act by multiplication), i.e. the $(K,H)$-bisets~$U$ for which there exists a biset morphism $U\to {_KG_H}$. \par
Indeed, a conjugation $(K,H)$-biset $U$ is a bifree $(K,H)$-biset isomorphic to a disjoint union of bisets of the form $(K\times H)/S$, where $S$ is a subgroup of $K\times H$ of the form
$$S_{g,A}=\{(^gx,x)\mid x\in A\}$$
where $A$ is a subgroup of $H$, and $g$ is an element of $G$ such that $^gA\leq K$. For such a transitive biset $(K\times H)/S$, the map
$$\forall (k,h)S\in (K\times H)/S,\; (k,h)S\mapsto kgh^{-1}$$
is a morphism of $(K,H)$-bisets. \par
Conversely, let $U$ be a $(K,H)$-biset for which there exists a biset morphism $\alpha:U\to {_KG_H}$. Then for any $u\in U$, the stabilizer $S_u$ of $u$ in $K\times H$ is the subgroup
$$S_u=\{(k,h)\in K\times H\mid k\cdot u\cdot h^{-1}=u\}$$
of $K\times H$. Then if $(k,h)\in S_u$,
$$\alpha(k\cdot u)=k\alpha(u)=\alpha(u\cdot h)=\alpha(u)h\mpoint$$
Let $A_u$ denote the projection of $S_u$ into $H$, and set $g_u=\alpha(u)$. It follows that $S_u\subseteq S_{g_u,A_u}$.\par
Conversely, if $(k,h)\in S_{g_u,A_u}$, then $k={^{g_u}}h$, and there exists some $x\in K$ such that $(x,h)\in S_u$, since $h\in A_u$. Thus $x\cdot u\cdot h^{-1}=u$, from which follows that
$$\alpha(x\cdot u)=xg_u=\alpha(u\cdot h)=g_uh\mvirg$$
hence $x={^{g_u}h}=k$, and $S_u=S_{g_u,A_u}$. Observation~\ref{observe} follows.\par
\npar In other words, conjugation $(K,H)$-bisets form a category ${\bf Conj}_{K,H}^G$, and there is a forgetful functor $\Phi:\bisetover{K}{H}{_KG_H}\to{\bf Conj}_{K,H}^G$ sending $(U,a)$ to $U$. This functor is full, preserves disjoint unions, and moreover it induces a surjection on the corresponding sets of isomorphism classes. This means that $\Phi$ induces a surjective group homomorphism (still denoted by $\Phi$) from $\mathcal{B}(_KG_H)$ to the Grothendieck group $\mathcal{B}_{K,H}^G$ of conjugation $(K,H)$-bisets.
\npar If $H$, $K$ and $L$ are subgroups of $G$, if $(U,a)$ is a $(K,H)$-biset over $_KG_H$ and $(V,b)$ is an $(L,K)$-biset over $_LG_K$, the composition $(V,b)\circ (U,a)$ is the $(L,H)$-biset over $_LG_H$ defined by the following diagram:
$$\xymatrix@R=1.5ex@C=2ex{
V\ar[dd]^-{b}&&U\ar[dd]^-{a}&&V\times_KU\ar[dd]^-{b\times_Ka}\\
&\circ&&=&\\
_LG_K&&_KG_H&&G\times_KG\ar[dd]^-{\mu}\\
&&&&\\
&&&&_LG_H
}$$
where $\mu$ is multiplication in $G$. This composition is associative, and additive with respect to disjoint unions. Hence it induces a composition
$$\widehat{\circ}:\mathcal{B}(_LG_K)\times\mathcal{B}(_KG_H)\to\mathcal{B}(_LG_H)\;\;.$$
Hence, one can define a category $\widehat{\mathbf{B}}(G)$ whose objects are the subgroups of~$G$, and such that ${\rm Hom}_{\widehat{\mathbf{B}}(G)}(H,K)=\mathcal{B}(_KG_H)$, for subgroups $H$ and $K$ of~$G$. Composition is given by $\widehat{\circ}$, and the identity morphism of the subgroup $H$ of $G$ in the category $\widehat{\mathbf{B}}(G)$ is the class of the biset $(_HH_H,i_H)$, where $i_H:{_HH_H}\to {_HG_H}$ is the inclusion map from $H$ to $G$.\par
Since the functor $\Phi$ maps the composition $\widehat{\circ}$ to the composition of bisets, and the identity morphism of $H$ in $\widehat{\mathbf{B}}(G)$ to the identity biset $_HH_H$, one can extend $\Phi$ to a functor $\widehat{\mathbf{B}}(G)\to\mathbf{B}(G)$, which is the identity on objects. \par
In other words, the category $\mathbf{B}(G)$ introduced in Section~3 of \cite{htw_Mackey-bisets} is the quotient of the category $\widehat{\mathbf{B}}(G)$ obtained by identifying morphisms which have the same image by $\Phi$.
\npar By the above Remark~\ref{ind}, when $H$ and $K$ are subgroups of $G$, there is a group isomorphism
$$\mathcal{B}(_KG_H)\cong \mathcal{B}\big(\Ind_{K\times H}^{G\times G}(_KG_H)\big)\;\;,$$
(with the usual identification of $(K,H)$-bisets with $(K\times H)$-sets). Now the biset $_KG_H$ is actually the restriction to $(K\times H)$ of the $(G,G)$-biset~$G$. By the Frobenius reciprocity, it follows that
$$\Ind_{K\times H}^{G\times G}(_KG_H)\cong \Ind_{K\times H}^{G\times G}\Res_{K\times H}^{G\times G}{(_GG_G)}\cong\big(\Ind_{K\times H}^{G\times G}\bullet\big)\times {_GG_G}\;\;,$$
where $\bullet$ is a set of cardinality 1. Since $\Ind_{K\times H}^{G\times G}\bullet\cong (G/K)\times(G/H)$, it follows (after switching $G/H$ and $G$) that
$$\Ind_{K\times H}^{G\times G}(_KG_H)\cong(G/K)\times G\times(G/H)\;\;,$$
where the $(G,G)$-biset structure of the right hand side is given by
$$\forall (a,b,x,y,g)\in G^5,\;\;a\cdot(xK,g,yH)\cdot b=(axK,agb,b^{-1}yH)\;\;.$$
\npar It should now be clear that the additive completion $\widehat{\mathbf{B}}_\bullet(G)$ is equivalent to the category whose objects are finite $G$-sets, where for any two finite $G$-sets $X$ and $Y$
$${\rm Hom}_{\mathbf{B}_\bullet(G)}(X,Y)=\mathcal{B}\big(_G(Y\times G\times X\big)_G)\;\;,$$
the $(G,G)$-biset structure on $(Y\times G\times X)$ being given as above by
$$\forall (a,b,g,x,y)\in G^3\times X\times Y,\;\;a\cdot(y,g,x)\cdot b=(ay,agb,b^{-1}x)\;\;.$$
Keeping track of the composition $\widehat{\circ}$ along the above isomorphism shows that the composition in the category $\widehat{\mathbf{B}}_\bullet(G)$ can be defined by linearity from the following: if $X$, $Y$, and $Z$ are finite $G$-sets, if 
$$\xymatrix@R=2ex{
&V\ar[ldd]_-{f}\ar[dd]^-{e}\ar[rdd]^-{d}&&&&U\ar[ldd]_-{c}\ar[dd]^-{b}\ar[rdd]^-{a}&\\
&&&\hbox{and}&&&\\
Z&G&Y&&Y&G&X
}
$$
are $(G,G)$-bisets over $(Z\times G\times Y)$ and $(Y\times G\times X)$, respectively, their composition is given by the following $(G,G)$-biset over $(Z\times G\times X)$
$$\xymatrix@R=4ex@C=2ex{
&(V\times_{d,c}U)/G\ar[ld]_-{\gamma}\ar[d]^-{\beta}\ar[rd]^-{\alpha}&\\
Z&G&X
}
$$
where $V\times_{d,c}U$ is the pullback of $V$ and $U$ over $Y$, i.e. the set of pairs $(v,u)\in V\times U$ with $d(v)=c(u)$, and $(V\times_{d,c}U)/G$ the set of orbits of $G$ on it for the action given by $(v,u)\cdot g=(vg,g^{-1}u)$. This makes sense because $d(v\cdot g)=g^{-1}d(v)=g^{-1}c(u)=c(g^{-1}\cdot u)$ if $d(v)=c(u)$. The map $(\gamma,\beta,\alpha)$ is given by
$$(\gamma,\beta,\alpha)\big((v,u)G\big)=\big(f(v),e(v)b(u),a(u)\big)\;\;.$$
\npar The functor $\Phi:\widehat{\mathbf{B}}(G)\to\mathbf{B}(G)$ extends uniquely to an additive functor $\Phi_\bullet:\widehat{\mathbf{B}}_\bullet(G)\to\mathbf{B}_\bullet(G)$, and the category $\mathbf{B}_\bullet(G)$ is the quotient of $\widehat{\mathbf{B}}_\bullet(G)$ obtained by identifying morphisms which have the same image by $\Phi_\bullet$. Clearly, two morphisms $f,g\in{\rm Hom}_{\widehat{\mathbf{B}}_\bullet(G)}(X,Y)$ are identified if and only if $f-g$ is in the kernel of the group homomorphism
$$\phi:\mathcal{B}\big({_G(Y\times G\times X)_G}\big)\to \mathcal{B}\big({_G(Y\times X)_G}\big)$$
induced by the correspondence
$$\xymatrix@R=2ex{
&U\ar[ldd]_-{c}\ar[dd]^-{b}\ar[rdd]^-{a}&&&&U\ar[ldd]_-{c}\ar[rdd]^-{a}&\\
&&&\mapsto &&&\\
Y&G&X&&Y&&X
}
$$
on bisets. In other words, a morphism $f$ in $\widehat{\mathbf{B}}_\bullet(G)$ gives the zero morphism in $\mathbf{B}_\bullet(G)$ if and only if it belongs to ${\rm Ker}\,\phi$. 
\npar Now the $(G,G)$-biset $_GG_G$ is isomorphic to $\Ind_{\Delta(G)}^{G\times G}\bullet$, where $\Delta(G)$ is the diagonal subgroup of $G\times G$. It follows that there is an isomorphism of $(G,G)$-bisets
$$Y\times G\times X\cong\Ind_{\Delta(G)}^{G\times G}(Y\times X)\;\;.$$
Hence, by Remark~\ref{ind}~again, since $\Delta(G)\cong G$,
$$\mathcal{B}\big({_G(Y\times G\times X)_G}\big)\cong\mathcal{B}\big({_G(Y\times X)}\big)\;\;,$$
where $_G(Y\times X)$ is the usual cartesian product with diagonal $G$-action. More precisely, this isomorphism is induced by the correspondence
$$\xymatrix@R=2ex{
&U\ar[ldd]_-{c}\ar[dd]^-{b}\ar[rdd]^-{a}&&&&\makebox[1ex]{$b^{-1}(1)$}\ar[ldd]_-{c}\ar[rdd]^-{a}&\\
&&&\mapsto&&&\\
Y&G&X&&Y&&X
}
$$
It is then easy to check that the composition of 
$$
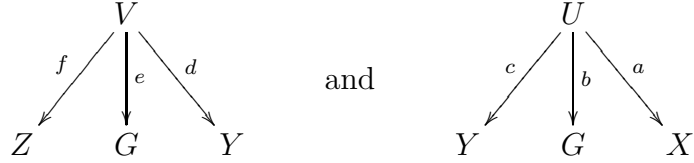
\xymatrix@R=2ex{
&V\ar[ldd]_-{f}\ar[dd]^-{e}\ar[rdd]^-{d}&&&&U\ar[ldd]_-{c}\ar[dd]^-{b}\ar[rdd]^-{a}&\\
&&&\hbox{and}&&&\\
Z&G&Y&&Y&G&X
}
$$
corresponds to the usual pullback diagram
$$
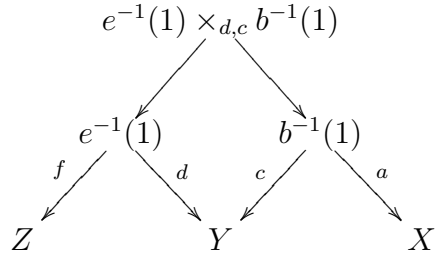
\xymatrix{
&&\makebox[1ex]{$e^{-1}(1)\times_{d,c}b^{-1}(1)\ar[ld]\ar[rd]$}&&\\
&\makebox[1ex]{$e^{-1}(1)$}\ar[ld]_-{f}\ar[rd]^-{d}&&\makebox[1ex]{$b^{-1}(1)$}\ar[ld]_-{c}\ar[rd]^-{a}&\\
Z&&Y&&X
}
$$
In other words, the category $\widehat{\mathbf{B}}_\bullet(G)$ is equivalent to the category $\mathbf{S}(G)$ whose objects are the finite $G$-sets, where 
$${\rm Hom}_{\mathbf{S}(G)}(X,Y)=\mathcal{B}\big({_G(Y\times X)}\big)\;\;,$$
and composition is induced by pullback. It has been shown by Lindner (\cite{lindner}, see also \cite{green}) that the additive functors on this category are precisely the Mackey functors for $G$.
\npar It remains to keep track of identifications by $\Phi$, i.e. to start with a morphism $f\in {\rm Hom}_{\mathbf{S}(G)}(X,Y)$, to lift it to 
$$f^+\in{\rm Hom}_{\widehat{\mathbf{B}}_\bullet(G)}(X,Y)=\mathcal{B}\big({_G(Y\times G\times X)_G}\big)\;\;,$$
and see when $f^+$ lies in ${\rm Ker}\,\phi$. Now $f$ is represented by a difference of two $G$-sets over $_G(Y\times X)$ of the form
$$
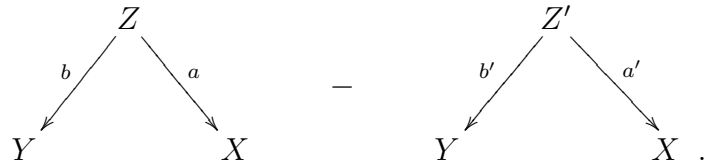
\xymatrix@R=2ex{
&Z\ar[ldd]_-{b}\ar[rdd]^-{a}&&&&Z'\ar[ldd]_-{b'}\ar[rdd]^-{a'}&\\
&&&-&&&\\
Y& &X&&Y& &X\;\;.
}
$$
By induction from $\Delta(G)$ to $G\times G$, the $G$-set on the left hand side lifts to the following  $(G\times G)$-set over $_{(G\times G)}(Y\times G\times X)$
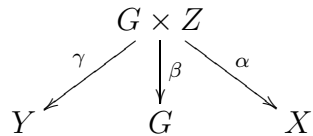
$$\xymatrix{
&G\times Z\ar[ld]_-{\gamma}\ar[d]^-{\beta}\ar[rd]^-{\alpha}&\\
Y&G&X
}
$$
where the $(G\times G)$-actions on $G\times Z$ and $Y\times G\times X$ are given respectively by $(s,t)\cdot (g,z)=(sgt^{-1},tz)$ and $(s,t)\cdot(y,g,x)=(sy,sgt^{-1},tx)$, and where
$$(\gamma,\beta,\alpha)(g,z)=\big(gb(z),g,a(z)\big)\;\;.$$
Similarly the $G$-set $\big(Z',(b',a')\big)$ lifts to $\big(G\times Z',(\gamma',\beta',\alpha')\big)$.\par
Now $f^+$ is in ${\rm Ker}\,\phi$ if and only if there is an isomorphism
$$\xymatrix@R=2ex{
&G\times Z\ar[ldd]_-{\gamma}\ar[rdd]^-{\alpha}&&&&G\times Z'\ar[ldd]_-{\gamma'}\ar[rdd]^-{\alpha'}&\\
&&&\stackrel{\theta}{\longrightarrow}&&&\\
Y& &X&&Y& &X\;\;.
}
$$
of $(G\times G)$-sets over $Y\times X$. Since $(g,z)=g\cdot(1,z)$ for any $(g,z)\in G\times Z$, it follows that $\theta$ is a map from $G\times Z$ to $G\times Z'$ of the form
$$(g,z)\mapsto \big(gu(z),v(z)\big)\;\;,$$
where $u$ is a map from $Z$ to $G$ and $v$ is a map from $Z$ to $Z'$. Now for any $(s,t)\in G\times G$, the equality
$$\theta\big((s,t)\cdot(g,z)\big)=(s,t)\cdot\theta\big((g,z)\big)$$
gives
$$\big(sgt^{-1}u(tz),v(tz\big)=\big(sgu(z)t^{-1},tv(z)\big)\;\;.$$
This is equivalent to
$$u(tz)={^tu(z)}\;\;\hbox{and}\;\;v(tz)=tv(z)\;\;.$$
This means that $u$ is a morphism of $G$-sets from $Z$ to $G^c$, which is the set $G$ with $G$-action by conjugation, and $v$ is a morphism of $G$-sets.\par
Moreover $\theta$ is a bijection if and only if $v$ is.\par
Finally $\theta$ is an morphism of $(G,G)$-bisets {\em over $Y\times X$} if and only if $\alpha'\circ \theta=a$ and $\gamma'\circ\theta=\gamma$, i.e. equivalently if
$$
a'\circ v=a\;\;\hbox{and}\;\;gu(z)\cdot b'\circ v(z)=g\cdot b(z)
$$
for any $(g,z)\in G\times Z$. In other words
$$a=a'\circ v\;\;\hbox{and}\;\;b=u*(b'\circ v)\;\;,$$
where, for any map $w:Z\to Y$, the map $u*w:Z\to Y$ is defined by
$(u*w)(z)=u(z)\cdot w(z)$. The map $u*w$ is a map of $G$-sets if $u:Z\to G^c$ and $w:Z\to Y$ are. Note that $w'=u*w$ if and only if $w=\bar{u}*w'$, where $\bar{u}:Z\to G^c$ is defined by $\bar{u}(z)=u(z)^{-1}$.\par
It follows that $f$ maps to the zero morphism in $\mathbf{B}(G)$ if and only if there exists $u:Z\to G^c$ and an isomorphism $v:Z\to Z'$ such that
$$a'\circ v=a\;\;\hbox{and}\;\;b'\circ v=u*b\;\;,$$
But then $v$ is an isomorphism
$$\xymatrix@R=2ex{
&Z\ar[ldd]_-{b'\circ v}\ar[rdd]^-{a'\circ v}&&&& Z'\ar[ldd]_-{b'}\ar[rdd]^-{a'}&\\
&&&\stackrel{v}{\longrightarrow}&&&\\
Y& &X&&Y& &X\;\;.
}
$$
of $G$-sets over $Y\times X$, and $f$ is also represented by the difference
$$\xymatrix@R=2ex{
&Z\ar[ldd]_-{b}\ar[rdd]^-{a}&&&&Z\ar[ldd]_-{u*b}\ar[rdd]^-{a}&\\
&&&-&&&\\
Y& &X&&Y& &X\;\;,
}
$$
since $a'\circ v=a$ and $b'\circ v=u*b$. These are the morphisms in the category $\mathbf{S}(G)$ that vanish in $\mathbf{B}_\bullet(G)$. In other words:
\begin{mth}{Theorem} \label{principal}Let $G$ be a finite group. Let $\sou{\mathbf{S}}(G)$ denote the quotient category of $\mathbf{S}(G)$ defined by setting, for any two finite $G$-sets $Y$ and $Y$
$$\Hom_{\sou{\mathbf{S}}(G)}(X,Y)=\mathcal{B}\big(_G(Y\times X)\big)/K(Y,X)\mvirg$$
where $K(Y,X)$ is the subgroup generated by the differences
\begin{equation}\label{diff}
\xymatrix@R=2ex{
&Z\ar[ldd]_-{b}\ar[rdd]^-{a}&&&&Z\ar[ldd]_-{u*b}\ar[rdd]^-{a}&\\
&&&-&&&\\
Y& &X&&Y& &X\;\;,
}
\end{equation}
where $a:Z\to X$, $b:Z\to Y$, and $u:Z\to G^c$ are morphisms of $G$-sets.\par
Then the functor $\Phi_\bullet$ induces an equivalence of categories $\sou{\mathbf{S}}(G)\cong \mathbf{B}_\bullet(G)$.
\end{mth}
Since the difference \ref{diff} factors as
$$\vcenter{\xymatrix@R=2ex@C=1.2ex{
&Z\ar[ldd]_-{b}\ar[rdd]^-{\rm Id}&\\
&&&\circ\\
Y&&Z&
}}
\left(\vcenter{
\xymatrix@R=2ex@C=1.2ex{
&Z\ar[ldd]_-{{\rm Id}}\ar[rdd]^-{\rm Id}&\\
&&&-\\
Z&&Z&
}}
\vcenter{
\xymatrix@R=2ex@C=1.2ex{
&Z\ar[ldd]_-{u*{\rm Id}}\ar[rdd]^-{\rm Id}&\!\!\!\\
&&&\rule{0ex}{1.3ex}\!\!\!\\
Z&&Z\!\!\!
}
}
\right)
\vcenter{\xymatrix@R=2ex@C=1.2ex{
&&Z\ar[ldd]_-{{\rm Id}}\ar[rdd]^-{a}&\\
\circ&&&&\\
&Z&&X&
}}
$$
the morphism vanishing in $\sou{\mathbf{S}}(G)$ are generated in the category $\mathbf{S}(G)$ by the morphisms of the form
$$\xymatrix@R=2ex{
&Z\ar[ldd]_-{\rm Id}\ar[rdd]^-{\rm Id}&&&&Z\ar[ldd]_-{u*\rm Id}\ar[rdd]^-{\rm Id}&\\
&&&-&&&\\
Z& &Z&&Z& &Z\;\;.
}
$$
\npar It follows that the additive functors from $\sou{\mathbf{S}}(G)$ to the category of abelian groups are exactly those Mackey functors (in the sense of Dress) such that for any $G$-set $Z$ and any $u:Z\to G^c$, the morphism $M_*(u*{\rm Id})$ is equal to the identity map of $M(Z)$. \par
This condition is additive with respect to $Z$, since the map $u*{\rm Id}_Z$ maps each $G$-orbit of $Z$ to itself. Hence these functors are exactly the functors for which the map $M_*(u*{\rm Id})$ is the identity map of $M(G/H)$, for any subgroup $H$ of $G$ and any $u:G/H\to G^c$. Such a map is of the form $gH\mapsto gcH$, where $c\in C_G(H)$. The map $u*{\rm Id}:G/H\to G/H$ is the map $gH\mapsto gcH$. \par
Translated in terms of the usual definition of Mackey functors, this map expresses the action of $c$ on $M(H)=M(G/H)$. This shows that additive functors from $\sou{\mathbf{S}}(G)$ to abelian groups are exactly the Mackey functors for the group $G$ such that, for any $H\leq G$, the centralizer $C_G(H)$ acts trivially on $M(H)$. These are the ``conjugation invariant Mackey functors" introduced in~\cite{htw_Mackey-bisets}.
\section{Fused $G$-sets}\label{fused G-sets}
Let $Z$ be any (finite) $G$-set. The multiplication $(u,v)\mapsto u*v$ endows the set ${\rm Hom}_{\gset{G}}(Z,G^c)$ with a group structure. Moreover, for any finite $G$-set~$X$, this group acts on the left on the set ${\rm Hom}_{\gset{G}}(Z,X)$, via $(u,f)\mapsto u*f$. This action is compatible with the composition of morphisms: if $Y$ is a finite $G$-set, if $u:Z\to G^c$ and $v:Y\to G^c$ are morphisms of $G$-sets, then for any morphisms of $G$-sets $f:Z\to Y$ and $g:Y\to X$, one checks easily that
\begin{equation}\label{composition}
(v*g)\circ(u*f)=\big(u*(v\circ f)\big)*(g\circ f)\mpoint
\end{equation}
\pagebreak[3]
\begin{mth}{Notation}
Let $\sgset{G}$ denote the category of {\em fused $G$-sets}: its objects are finite $G$-sets, and for any finite $G$-sets $Z$ and $Y$
$${\rm Hom}_{\isgset{G}}(Z,Y)={\rm Hom}_{\igset{G}}(Z,G^c)\dom{\rm Hom}_{\igset{G}}(Z,Y)\;\;.$$
The composition of morphisms in $\sgset{G}$ is induced by the composition of morphisms in $\gset{G}$.
\end{mth}
\begin{rem}{Remark} For any $G$-set $Y$, set $Y^I=Y\times G^c$. This notation is chosen to evoke a path object in homotopy theory (cf. \cite{dwyer-spalinski} Section 4.12). There is a natural morphism $p:Y^I\to Y\times Y$, defined by $p(y,g)=(y,gy)$, for $y\in Y$ and $g\in G$, and a morphism $i:Y\to Y^I$ defined by $i(y)=(y,1)$, for $y\in Y$. The composition $p\circ i$ is equal to the diagonal map $Y\to Y\times Y$.\par
Two morphisms $a,b:Z\to Y$ in $\gset{G}$ are equal in the category $\sgset{G}$ if and only if the morphism $(a,b):Z\to Y\times Y$ factors as
$$\xymatrix{
&Y^I\ar[d]^-p\\
Z\ar[ur]^-\varphi\ar[r]_-{(a,b)}&Y\times Y
}
$$
for some morphism of $G$-sets $\varphi:Z\to Y^I$.
\end{rem}
\begin{rem}{Remark} It follows from \ref{composition} that the map $u\mapsto u*\Id_Z$ is a group antihomomorphism from $\Hom_{\igset{G}}(Z,G^c)$ to the group of $G$-automorphisms of~$Z$. Hence a morphism $\sou{f}:Z\to Y$ in the category $\sgset{G}$ is an isomorphism if and only if any of its representatives $f:Z\to Y$ in $\gset{G}$ is an isomorphism.
\end{rem}
\masubsect{Weak pullbacks of fused $G$-sets} Disjoint union of $G$-sets is a coproduct in $\sgset{G}$. There is also a weak version of pullback in $\sgset{G}$~: let
$$\xymatrix{
&T\ar[dl]_-{\sou{c}}\ar[dr]^-{\sou{d}}&\\
X\ar[dr]_-{\sou{a}}&&Y\ar[dl]^-{\sou{b}}\\
&Z&
}$$
be a commutative diagram in $\sgset{G}$, where underlines denote the images in $\sgset{G}$ of morphisms in $\gset{G}$. This means that $\sou{a}\circ\sou{c}=\sou{b}\circ\sou{d}$, i.e. that there exists $u\in\Hom_{\igset{G}}(T,G^c)$ such that
$$b\circ d=u*(a\circ c)\mpoint$$
But $u*(a\circ c)=a\circ(u*c)$. It follows that there is a unique morphism $e\in\Hom_{\igset{G}}(T,X\times_{a,b}Y)$ such that the diagram
$$
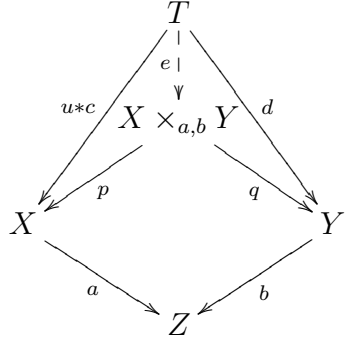
\xymatrix{
&T\ar[ddl]_-{u*c}\ar@{-->}[d]_-{e}\ar[ddr]^-{d}\\
&X\times_{a,b}Y\ar[dl]^-{{p}}\ar[dr]_-{{q}}&\\
X\ar[dr]_-{{a}}&&Y\ar[dl]^-{{b}}\\
&Z&
}$$
is commutative in $\gset{G}$, where $p:X\times_{a,b}Y\to X$ and $q:X\times_{a,b}Y\to Y$ are the canonical morphisms from the pullback $X\times_{a,b}Y$. In other words, the diagram
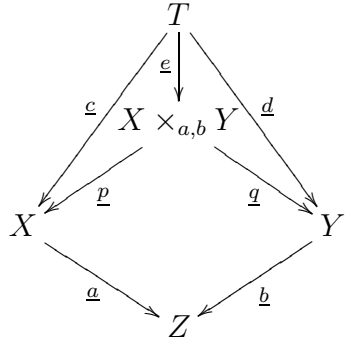
\begin{equation}\label{diagram}\xymatrix{
&T\ar[ddl]_-{\sou{c}}\ar[d]_-{\sou{e}}\ar[ddr]^-{\sou{d}}\\
&X\times_{a,b}Y\ar[dl]^-{\sou{p}}\ar[dr]_-{\sou{q}}&\\
X\ar[dr]_-{\sou{a}}&&Y\ar[dl]^-{\sou{b}}\\
&Z&
}
\end{equation}
is commutative in $\sgset{G}$.\par
But still $(X\times_{a,b}Y,\sou{p},\sou{q})$ need not be a pullback in $\sgset{G}$, since the morphism~$\sou{e}$ making Diagram~\ref{diagram} commutative is generally not unique, as $e$ itself depends on the choice of $u$. Moreover, the lifts $a$ and $b$ of $\sou{a}$ and $\sou{b}$ to $\gset{G}$ are not unique~: it should be noted however that if $a'=v*a$ and $b'=w*b$ are other lifts of $a$ and $b$, respectively, where $v\in\Hom_{\igset{G}}(X,G^c)$ and $w\in\Hom_{\igset{G}}(Y,G^c)$, then the map $f:(x,y)\mapsto \big(v(x)x,w(y)y\big)$ is an isomorphism of $G$-sets from $X\times_{a',b'}Y$ to $X\times_{a,b}Y$, such that the diagram
\newlength{\prov}
\settowidth{\prov}{$X\otimes_{a,b}Y$}
$$\xymatrix{
& \makebox[\prov]{$X\times_{a',b'}Y$}\ar[dl]_-{p'}\ar[dr]^<(0.6){q'}\ar[drrr]^-{f}&&&&\\
X\ar[dr]_-{a'}\ar[drrr]_<(0.7){v*\Id}&&Y\ar[dl]|<(0.48)\hole^<(0.1){b'}\ar[drrr]|<(0.43)\hole_<(0.6){w*\Id}&&X\times_{a,b}Y\ar[dl]_<(0.7){p}\ar[dr]^-{q}&\\
&Z\ar[drrr]_-{\Id}&&X\ar[dr]_<(0.3){a}&&Y\ar[dl]^-{b}\\
&&&&Z&
}
$$
is commutative in $\gset{G}$. Since $\sou{a'}=\sou{a}$, $\sou{b'}=\sou{b}$, $\sou{v*\Id}=\sou{\Id}$, and $\sou{w*\Id}=\sou{\Id}$, this yields a commutative diagram
$$\xymatrix{
&X\times_{a',b'}Y\ar[rr]^-{\sou{f}}\ar[dl]_-{\sou{p'}}\ar[drrr]|<(0.2)\hole_{\sou{q'}}&&X\times_{a,b}Y\ar[dlll]^{\sou{p}}\ar[dr]^{\sou{q}}&\\
X\ar[drr]_{\sou{a}}&&&&Y\ar[dll]^{\sou{b}}\\
&&Z&&
}
$$
in $\sgset{G}$, and $\sou{f}$ is an isomorphism. This shows that the weak pullback $X\times_{a,b}Y$ only depends on $\sou{a}$ and $\sou{b}$ in the category $\sgset{G}$. For this reason, it may be denoted by $X\times_{\sou{a},\sou{b}}Y$.
\masubsect{Spans of fused $G$-sets}
Recall (cf. \cite{yoneda}, \cite{benabou} for the general definition) that if $X$ and $Y$ are finite $G$-sets, then a {\em span} $\Lambda_{Z,\sou{a},\sou{b}}$ over $X$ and $Y$ in the category $\sgset{G}$ is a diagram of the form
$$
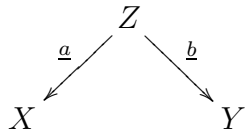
\xymatrix{
&Z\ar[dl]_-{\sou{a}}\ar[dr]^-{\sou{b}}&\\
X&&Y
}
$$
where $Z$ is a finite $G$-set and $\sou{a}, \sou{b}$ are morphisms in the category $\sgset{G}$. Two spans $\Lambda_{Z,\sou{a},\sou{b}}$ and $\Lambda_{Z',\sou{a'},\sou{b'}}$ over $X$ and $Y$ are equivalent if there exists an isomorphism $\sou{f}:Z\to Z'$ in $\sgset{G}$ such that the diagram
$$
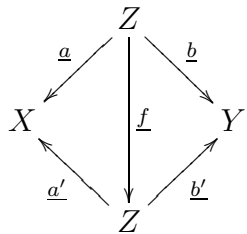
\xymatrix{
&Z\ar[dd]^-{\sou{f}}\ar[dl]_-{\sou{a}}\ar[dr]^-{\sou{b}}&\\
X&&Y\\
&Z\ar[ul]^-{\sou{a'}}\ar[ur]_-{\sou{b'}}&
}
$$
is commutative. The set of equivalence classes of spans of fused $G$-sets over $X$ and $Y$ is an additive monoid, where the addition is defined by disjoint union (i.e. $\Lambda_{Z_1,\sou{a}_1,\sou{b}_1}+\Lambda_{Z_2,\sou{a}_2,\sou{b}_2}=\Lambda_{Z_1\sqcup Z_2,\sou{a}_1\sqcup \sou{a}_2,\sou{b}_1\sqcup \sou{b}_2}$). The corresponding Grothendieck group is isomorphic to $\Hom_{\sou{\mathbf{S}}(G)}(Y,X)$. \par
It should be noted that even if there is no pullback construction in the category $\sgset{G}$, the isomorphism classes of spans in $\sgset{G}$ can still be composed by {\em weak pullback}, and this induces the composition of morphisms in~$\sou{\mathbf{S}}(G)$.\par
\section{Fused Mackey functors}\label{fused Mackey functors}
\begin{mth}{Definition} \label{fused Mackey}Let $R$ be a commutative ring. Let $R\,\mathbf{S}(G)$ (resp. $R\,\sou{\mathbf{S}}(G)$) denote the $R$-linear extension of the category $\mathbf{S}(G)$ (resp. $\sou{\mathbf{S}}(G)$), defined as follows:
\begin{itemize}
\item The objects of $R\,\mathbf{S}(G)$ and $R\,\sou{\mathbf{S}}(G)$ are finite $G$-sets.
\item For finite $G$ sets $X$ and $Y$, 
$$\Hom_{R\,{\mathbf{S}}(G)}(X,Y)=R\otimes_\Z\Hom_{{\mathbf{S}}(G)}(X,Y)\mvirg$$
$$\Hom_{R\,\sou{\mathbf{S}}(G)}(X,Y)=R\otimes_\Z\Hom_{\sou{\mathbf{S}}(G)}(X,Y)\mpoint$$
\item Composition of morphisms is induced by the pullback in $\gset{G}$ (resp. the weak pullback in $\sgset{G}$).
\end{itemize}
A Mackey functor for $G$ over $R$ in the sense of Lindner (\cite{lindner}) is an $R$-linear functor from $R\,\mathbf{S}(G)$ to the category $\gMod{R}$ of $R$-modules.\par
Similarly, a {\em fused Mackey functor} for $G$ over $R$ is an $R$-linear functor from $R\,\sou{\mathbf{S}}(G)$ to $\gMod{R}$. A morphism of fused Mackey functors is a natural transformation of functors. Fused Mackey functors for $G$ over $R$ form a category denoted by $\mathsf{Mack}_R^f(G)$.
\end{mth}
The following is an equivalent definition of fused Mackey functors, {\em \`a la } Dress: 
\begin{mth}{Definition} Let $R$ be a commutative ring. A {\em fused Mackey functor} for the group $G$ over $R$ is a bivariant $R$-linear functor $M=(M^*,M_*)$ from $\sgset{G}$ to $\gMod{R}$ such that:
\begin{enumerate}
\item For any finite $G$-sets $X$ and $Y$, the maps
$$\xymatrix{
M(X)\oplus M(Y)\ar[rrr]<1ex>^-{(M_*(\sou{i}_X),M_*(\sou{i}_Y)}&&&M(X\sqcup Y)\ar[lll]<1ex>^-{(M^*(\sou{i}_X),M^*(\sou{i}_Y)}
}
$$ 
induced by the canonical inclusions $i_X:X\to X\sqcup Y$ and $i_Y:Y\to X\sqcup Y$ are mutual inverse isomorphisms.
\item If 
$$\xymatrix{
&X\times_{\sou{a},\sou{b}}Y\ar[dl]^-{\sou{p}}\ar[dr]_-{\sou{q}}&\\
X\ar[dr]_-{\sou{a}}&&Y\ar[dl]^-{\sou{b}}\\
&Z&
}
$$
is a weak pullback diagram in $\sgset{G}$, then $M^*(\sou{a})M_*(\sou{b})=M_*(\sou{p})M^*(\sou{q})$.
\end{enumerate}
A morphism of fused Mackey functors is a natural transformation of bivariant functors.
\end{mth}
The category $\mathsf{Mack}_R^f(G)$ can be viewed as a full subcategory of the category $\mathsf{Mack}_R(G)$ of Mackey functors for~$G$ over $R$. In the case $R=\Z$, this category is equivalent to the category of conjugation invariant Mackey functors introduced in~\cite{htw_Mackey-bisets}.\par
The inclusion functor $\mathsf{Mack}_R^f(G)\hookrightarrow \mathsf{Mack}_R(G)$ has a left adjoint:
\begin{mth}{Definition} \label{adjoint}Let $M$ be a Mackey functor for $G$ over $R$, in the sense of Lindner, i.e. an $R$-linear functor $R\,\mathbf{S}(G)\to \gMod{R}$. When $X$ is a finite $G$-set, set
$$M^f(X)=M(X)/\sum_{Z,a,u} \Im \big(M(\Lambda_{a,\Id_Z})- M(\Lambda_{u*a,\Id_Z})\big)\mvirg$$
where the summation runs through triples $(Z,a,u)$ consisting of a finite $G$-set~$Z$, and morphisms of $G$-sets $a:Z\to X$ and $u:Z\to G^c$, and $\Lambda_{a,\Id_Z}$ denotes the span
$$\xymatrix@C=2.5ex{
&Z\ar[dl]_-a\ar@{=}[dr]^-{\Id_Z}&\\
X&&Z
}
$$
of $G$-sets.
\end{mth}
\begin{mth}{Proposition} Let $R$ be a commutative ring, and $G$ be a finite group.
\begin{enumerate}
\item Let $M$ be a Mackey functor for $G$ over $R$. The correspondence 
$$X\mapsto M^f(X)$$
is a fused functor $M^f$ for $G$ over $R$.
\item The correspondence $\mathcal{F}:M\mapsto M^f$ is a functor from $\mathsf{Mack}_R(G)$ to $\mathsf{Mack}_R^f(G)$, which is left adjoint to the inclusion functor  
$$\mathcal{I}:\mathsf{Mack}_R^f(G)\hookrightarrow \mathsf{Mack}_R(G)\mpoint$$
Moreover $\mathcal{F}\circ\mathcal{I}$ is isomorphic to the identity functor of $\mathsf{Mack}_R^f(G)$. 
\end{enumerate}
\end{mth}
\pf For Assertion 1, to prove that $M^f$ is a Mackey functor, observe that if $\Lambda_{Z,a,b}$ is a span of finite $G$-sets of the form
$$\xymatrix@C=2.5ex{
&Z\ar[dl]_-a\ar[dr]^-{b}&\\
X&&Y
}
$$
and $u:Z\to G^c$ is a morphism of $G$-sets, then
$$\Lambda_{Z,a,b}-\Lambda_{Z,u*a,b}=(\Lambda_{Z,a,\Id_Z}-\Lambda_{Z,u*a,\Id_Z})\circ\Lambda_{Z,\Id_Z,b}\mpoint$$
It follows that the $R$-module
$$\sum_{Z,a,u} \Im \big(M(\Lambda_{a,\Id_Z})- M(\Lambda_{u*a,\Id_Z})\big)$$
is equal to the sum
$$\sum_{Z,a,b,u} \Im \big(M(\Lambda_{a,b})- M(\Lambda_{u*a,b})\big)\mpoint$$
In other words, it is equal to the image by $M$ of the $R$-submodule $K_R(X,Y)$ of $\Hom_{R\,\mathbf{S}(G)}(Y,X)$ generated by the morphisms $\Lambda_{a,b}-\Lambda_{u*a,b}$, i.e. to the kernel of the quotient morphism
$$\Hom_{R\,\mathbf{S}(G)}(Y,X)\to \Hom_{R\,\sou{\mathbf{S}}(G)}(Y,X)\mpoint$$
This shows that $K_R$ is an ideal in the category $R\,\mathbf{S}(G)$. So if $M$ is an $R$-linear functor $R\,\mathbf{S}(G)\to\gMod{R}$, the correspondence 
$$X\mapsto M^f(X)=M(X)/\sum_{f\in K_R(X,Y)}\Im M(f)$$ 
is an $R$-linear functor from the quotient category $R\,\sou{\mathbf{S}}(G)$ to $\gMod{R}$.\par
Assertion 2 is straightforward: first it is clear that $\mathcal{F}\circ\mathcal{I}$ is isomorphic to the identity functor, since $N^f=N$ when $N$ is a fused Mackey functor. This isomorphism $\mathcal{F}\circ\mathcal{I}\cong \Id_{\mathsf{Mack}_R^f(G)}$ provides the counit of the adjunction. Next for any Mackey functor $M$, there is a projection morphism $M\to \mathcal{I}\mathcal{F}(M)$, and this yields the unit of the adjunction.\findemo
\begin{rem}{Remark} Assertion 2 shows that $\mathsf{Mack}_R^f(G)$ is a {\em reflective} subcategory of $\mathsf{Mack}_R^f(G)$ (cf. \cite{cwm}, Chapter \romain{4}, Section{3}).
\end{rem}
\begin{rem}{Remark} If the Mackey functor $M$ is given in the sense of Dress, then for any finite $G$-set $X$
$$M^f(X)=M(X)/\sumb{a:Z\to X}{u:Z\to G^c}\Im \big(M_*(a)-M_*(u*a)\big)\mvirg$$
where $Z$ is a finite $G$-set, and $a,u$ are morphisms of $G$-sets.
\end{rem}
\pagebreak[3]
\begin{mth}{Corollary} \label{enough proj}\begin{enumerate}
\item If $P$ is a projective Mackey functor, then $P^f$ is projective in the category $\mathsf{Mack}_R^f(G)$.
\item The category $\mathsf{Mack}_R^f(G)$ has enough projective objects. More precisely, if $N$ is a fused Mackey functor, and $\theta:P\to \mathcal{I}(N)$ is an epimorphism in $\mathsf{Mack}_R(G)$ from a projective Mackey functor $P$, then $\mathcal{F}(\theta):P^f\to N$ is an epimorphism in $\mathsf{Mack}_R^f(G)$.
\end{enumerate}
\end{mth}
\pf Assertion 1 follows from the fact that $\mathcal{F}$ is left adjoint to the exact functor $\mathcal{I}$. Assertion 2 is then straightforward.\findemo
\section{The fused Mackey algebra}\label{algebra}
When $G$ is a finite group, set $\Omega_G=\mathop{\sqcup}_{H\leq G}\limits G/H$, and let $RB_{\Omega_G}$ denote the Dress construction for the Burnside functor $RB$ over the ring $R$. Recall that $RB_{\Omega_G}$, as a Mackey functor in the sense of Dress, is obtained by precomposition of $RB$ with the endofunctor $X\mapsto X\times\Omega_G$ of $\gset{G}$.\par
Also recall (cf. \cite{green} Lemma 7.3.2 and Proposition 4.5.1) that the functor $RB_{\Omega_G}$ is a progenerator of the category $\mathsf{Mack}_R(G)$, and that the algebra $\End_{\mathsf{Mack}_R(G)}(B_{\Omega_G})\cong B(\Omega_G^2)$ is isomorphic to the Mackey algebra $\mu_R(G)$ of $G$ over $R$, introduced by Th\'evenaz and Webb (\cite{thevwebb}).\par
It follows from Corollary~\ref{enough proj} that the functor $(RB_{\Omega_G})^f$ is a progenerator in the category $\mathsf{Mack}_R^f(G)$. Hence this category is equivalent to the category of modules over the algebra $\End_{\mathsf{Mack}_R^f(G)}\big((RB_{\Omega_G})^f\big)$.
\begin{mth}{Definition} The {\em fused Mackey algebra}  of $G$ over $R$ is the algebra
$$\mu_R^f(G)=\End_{\mathsf{Mack}_R^f(G)}\big((RB_{\Omega_G})^f\big)\mpoint$$
\end{mth}
\begin{mth}{Lemma} \label{yoneda}Let $X$ be a finite $G$-set. Then $(RB_X)^f$ is isomorphic to the Yoneda functor $\Hom_{R\,\sou{\mathbf{S}}(G)}(X,-)$.
\end{mth}
\pf Denote by $\mathcal{Y}_X$ the Yoneda functor $\Hom_{R\,\sou{\mathbf{S}}(G)}(X,-)$. For any fused Mackey functor $N$ for $G$ over $R$
\begin{eqnarray*}
\Hom_{\mathsf{Mack}_R^f(G)}\big((RB_X)^f,N)&\cong& \Hom_{\mathsf{Mack}_R(G)}\big(RB_X,\mathcal{I}(N)\big)\\
&\cong& \mathcal{I}(N)(X)\cong N(X)\\
&\cong&\Hom_{\mathsf{Mack}_R^f(G)}(\mathcal{Y}_X,N)\mpoint
\end{eqnarray*}
The lemma follows, since all these isomorphisms are natural.\findemo
\begin{mth}{Theorem} The fused Mackey algebra $\mu_R^f(G)$ is isomorphic to the quotient of the algebra $RB(\Omega_G^2)\cong\mu_R(G)$ by the $R$-module generated by differences of the form 
$$
\xymatrix@C=2.5ex@R=2ex{
&Z\ar[ldd]_-{b}\ar[rdd]^-{a}&&&&Z\ar[ldd]_-{u*b}\ar[rdd]^-{a}&\\
&&&-&&&\\
\Omega_G& &\Omega_G&&\Omega_G& &\;\Omega_G\;,
}
$$
where $a,b: Z\to\Omega_G$ and $u:Z\to G^c$ are morphisms of $G$-sets.
\end{mth}
\pf This follows from Lemma~\ref{yoneda}, since the quotient in the theorem is precisely $\End_{R\,\sou{\mathbf{S}}(G)}(\Omega_G)$.\findemo 
\begin{rem}{Remark} One can deduce from this theorem that the fused Mackey algebra $\mu_R^f(G)$ is always free of finite rank as an $R$-module, and this rank does not depend on the commutative ring $R$. More precisely, Th\'evenaz and Webb have shown (\cite{thevwebb} Proposition 3.2) that the Mackey algebra $\mu_R(G)$ has an $R$-basis consisting of elements of the form
$$t_K^H\,c_{g,K}\,r_{K^g}^L\mvirg$$
where $(H,L,g,K)$ runs through a set of representatives of 4-tuples consisting of two subgroups $H$ and $L$ of $G$, and element $g$ of $G$, and a subgroup $K$ of $H\cap{^gL}$, for the equivalence relation $\equiv$ given by
$$(H,L,g,K)\equiv (H',L',g',K') \Leftrightarrow \left\{\begin{array}{l}H=H',\;L=L',\\
\hbox{and}\\
\exists h\in H,\;\exists l\in L,\;g'=hgl,\;K'={^hK}\mpoint\end{array}\right.$$
Similarly, the quotient algebra $\mu_R^f(G)$ of $\mu_R(G)$ has a basis consisting of the images of the elements $t_K^H\,c_{g,K}\,r_{K^g}^L$, where $(H,L,g,K)$ runs through a set of representatives of 4-tuples as above, modulo the relation $\equiv^f$ defined by
$$(H,L,g,K)\equiv^f (H',L',g',K') \Leftrightarrow \left\{\begin{array}{l}H=H',\;L=L',\\
\hbox{and}\\
\exists h\in H,\;\exists l\in L,\:\exists x\in C_G(K),\\
\;\;g'=hxgl,\;K'={^hK}\mpoint\end{array}\right.$$

\end{rem}

\centerline{\rule{5ex}{.1ex}}
\begin{flushleft}
Serge Bouc - CNRS-LAMFA, Universit\'e de Picardie, 33 rue St Leu, 80039, Amiens Cedex 01 - France. \\
{\tt email : serge.bouc@u-picardie.fr}\\
{\tt web~~ : http://www.lamfa.u-picardie.fr/bouc/}
\end{flushleft}

\end{document}